\documentclass[12pt]{article}
\usepackage[T2A,T1]{fontenc}
\usepackage{textcomp}
\usepackage{amsthm}
\usepackage{amsmath}
\usepackage{amssymb}
\usepackage{esint}
\usepackage{xcolor}
\newcommand {\RR} {\mathbb R}
\newcommand {\TT} {\mathbb T}
\newcommand {\ZZ} {\mathbb Z}
\newcommand {\CC} {\mathbb C}
\newcommand {\NN} {\mathbb N}
\usepackage[english,russian]{babel}
\makeatletter

\theoremstyle{plain}
\newtheorem{thm}{\protect\theoremname}
\newtheorem{lemma}{Lemma}
\ifx\proof\undefined
\newenvironment{proof}[1][\protect\proofname]{\par
\normalfont\topsep6\p@\@plus6\p@\relax
\trivlist
\itemindent\parindent
\item[\hskip\labelsep
\scshape
#1]\ignorespaces
}{%
\endtrivlist\@endpefalse
}
\providecommand{\proofname}{Proof}
\fi

\makeatother

\usepackage{babel}
\providecommand{\theoremname}{Theorem}

\begin{document}

\title{Stechkin's problem for functions of a self-adjoint operator in a Hilbert space, Taikov-type inequalities and their applications
\thanks{This project was supported by Simons Collaboration Grant N. 210363}
}


\author{Vladyslav Babenko$^1$,
        Yuliya Babenko$^2$, Nadiia Kriachko$^3$
}


\date{$^{1,3}$Department of Mathematics and Mechanics, Dnepropetrovsk National University, Gagarina pr., 72, Dnepropetrovsk, 49010, Ukraine \\
           $^2$Department of Mathematics, Kennesaw State University, 1100 South Marietta Pkwy, MD \# 9085, Marietta, GA, 30060, USA}

\maketitle

\begin{abstract}
 In this paper we solve the problem of approximating functionals $(\varphi(A)x, f)$ (where $\varphi(A)$ is some function of self-adjoint operator $A$) on the class of elements of a Hilbert space that is defined with the help of another function $\psi (A)$ of the operator $A$. In addition, we obtain a series of sharp Taikov-type additive inequalities that estimate $|(\varphi(A)x, f)|$ with the help of  $\| \psi (A)x\|$ and $\| x\|$. We also present several applications of the obtained results. First, we find sharp constants in inequalities of the type used in H${\rm{\ddot{o}}}$rmander theorem on comparison of operators in the case when operators are acting in a Hilbert space and are functions of a self-adjoint operator. As another application we obtain Taikov-type inequalities for functions of the operator $\frac1i \frac {d}{dt}$ in the spaces $L_2(\RR)$ and $L_2(\TT)$, as well as for integrals with respect to spectral
measures, defined with the help of classical orthogonal polynomials.\newline
{\bf Keywords}{functions of operators \and Landau-Kolmogorov inequalities \and best approximation \and orthogonal polynomials}\newline
{\bf Mathematics Subject Classification} (2000) MSC 26D10,  MSC 47A63,
MSC 41A17, MSC 47A58
\end{abstract}

\section{Introduction.}

We begin by stating Stechkin's problem (see, \cite{Stechkin_2}) of best approximation of an operator by linear bounded operators on some class of elements.

Let $X, \, Y$ be Banach spaces; $A:X\to Y$ be some operator (not necessarily linear) with a domain $D(A)\subset X$. Let also ${\cal L}(N)={\cal L}(N;X,Y)$ be a set of linear bounded operators $T:X\to Y$, with norms bounded above by $N>0$, and $Q\subset D(A)$ be some class of elements. The quantity
\[
U(T)=U(A,Q,T):=\sup\{ \| Ax-Tx\|_Y\; :\; x\in Q\}
\]
is called the {\it deviation} of the operator $T\in {\cal L}(N)$ from operator $A$ on the class $Q$, and
\begin{equation}\label{appr_oper}
E(N)=E(A,Q,N):=\inf\{ U(T)\; :\; T\in {\cal L}(N)\}
\end{equation}
is called the {\it best approximation} of the operator $A$ by the set of operators ${\cal L}(N)$ on the class $Q$.

The problem is to compute the quantity $E(N)=E(A,Q,N)$ and find an operator $T\in {\cal L}(N)$ that realizes the infimum on the right-hand side of (\ref{appr_oper}).

This problem has actively being studied by many authors. The history of the question, survey of the known results as well as further references can be found in \cite{Arestov}, \cite[Ch. 7]{babenko_inequalities}.
The high interest in this problem has been partially due to its close connection to Landau-Kolmogorov-type inequalities, which provide an estimate for the norm of an intermediate derivative of a function with the help of norms of the function itself and norms of higher derivatives. Such inequalities, as well as their various generalizations, found applications in various areas (see, for instance, \cite{Arestov}, \cite{babenko_inequalities}).

We demonstrate this connection on additive inequalities of Landau-Kolmogorov type. In addition, we present here one of the methods to solve Stechkin's problem.

We assume that operator $A$ is homogeneous. In addition, let a linear operator $B$ from a space $X$ into a Banach space $Z$ such that $D(B)\subset D(A)$ be given.
Usually, the class $Q$ in Stechkin's problem is taken to be
\[
Q=W^B=\{ x\in D(B)\; :\; \| Bx\|_Z\le 1\}.
\]

For any $N>0$, arbitrary operator $T\in {\cal L}(N)$, and arbitrary $x\in D(B)$ the following inequality holds
\begin{equation}\label{equality}
\| Ax\|_Y\le U(T)\| Bx\|_Z+N\| x\|_X,
\end{equation}
which implies
\begin{equation}\label{equality2}
\| Ax\|_Y\le E(N)\| Bx\|_Z+N\| x\|_X.
\end{equation}
The last inequality is a Landau-Kolmogorov{-type} inequality in additive form.
If for some operator $T$ and some element $\overline{x}\in D(B)$ inequality (\ref{equality}) becomes an equality, then comparing it with inequality (\ref{equality2}) written for $\overline{x}$, we see that $E(N)\ge U(T)$, and hence
\[
E(N)= U(T).
\]
Therefore, operator $T$ delivers the infimum in the right-hand side of (\ref{appr_oper}).

Let $G$ be the real line $\mathbb{R}$ or unit circle $\mathbb{T}.$
By $L_{p}(G)$ ($1\le p\le \infty$) we denote the space of all functions $x\,:\, G\rightarrow\mathbb{C}$,
such that
$$
\left\Vert x\right\Vert _{L_{p}(G)}:=
\left\{\begin{array}{clcr}
\mbox{$\left(\displaystyle \int\limits_{G}\left| x(t)\right|^{p}dt\right)^{\frac{1}{p}}$}, & \mbox{ если $1\leq p<\infty,$} \\
\mbox{${\rm esssup}\left\{ \left| x(t)\right| :\, t\in G\right\}$} , & \mbox{ если $p=\infty$}.
\end{array}\right.
$$
By $L_0(G)$ we denote the space of all measurable and almost everywhere finite functions $x\; :\; G\to \CC$. By $C(G)$ we denote the space of all continuous bounded functions $x\,:\, G\:\rightarrow\:\mathbb{C}$ equipped with the uniform norm.

By $L_{p,s}^{r}(G),\, r\in\mathbb{N}, \; 1\le p,s\le\infty,$ we denote the space of all functions $x\in L_{p}(G),$ such that their $(r-1)$-{\it st} derivative is locally absolutely continuous and $r$-{\it th} derivative is in $L_{s}(G).$ {Set} $W^{r}_{p,s}=\{ x\in L_{p,s}^{r}(G)\; :\; \| x^{(r)}\|_{L_{s}(G)}\le 1\}$.

In 1968 Taikov \cite{Taikov} considered Stechkin's problem in the following setting $X=Z=L_2(\RR),\; Y=C(\RR)$, $A=\frac {d^k}{dt^k}$, $B=\frac {d^r}{dt^r}$ (so that $W^B=W^{r}_{2,2}$),
$k,\, r\in\mathbb{N},\; k<r$. Taikov proved that if
\[
a=\left\{ \frac{r-k-1/2}{2r^{2}}\frac{1}{\sin\pi\frac{2k+1}{2r}}\right\} ^{1/2},\;\; b=\left\{ \frac{r+1/2}{2r^{2}}\frac{1}{\sin\pi\frac{2k+1}{2r}}\right\} ^{1/2}
\]
and $N=ah^{-k-1/2},\: h>0,$ then
\[
E\left(\frac{d^k}{dt^r},W_{2,2}^r,N\right)=bh^{r-k-1/2}.
\]
Clearly, the problem of computing the quantity $E\left( \frac {d^k}{dt^k},W^{r}_{2,2},N\right)$ is equivalent to the problem of approximating the unbounded functional $x^{(k)}(0)$ on the class $W_{2,2}^{r}(\mathbb{R})$
by bounded functionals (case of $Y=\RR$). In addition, in \cite{Taikov} Taikov obtained best {possible} inequalities estimating the uniform norm of k$th$, $0<k<r,$ derivative
of the function $x\in L_{2,2}^{r}(\mathbb{R})$ with the help of $L_{2}(\RR)$-norm{s} of $x$ and $x^{(r)}$ in additive and multiplicative forms. Let us present the additive inequality here. For any $h>0$
\begin{equation}
\left\Vert x^{(k)}\right\Vert _{C(\mathbb{R})}\leq ah^{-k-1/2}\left\Vert x\right\Vert _{L_{2}(\mathbb{R})}+bh^{r-k-1/2}\left\Vert x^{(r)}\right\Vert _{L_{2}(\mathbb{R})}.\label{eq:2}
\end{equation}

In 2012, Babenko and Bilichenko \cite{Babenko_Bilichenko} solved the problem of best approximation of the functional $F_{f}(x):=\left(A^{k}x,\, f\right)$
on the class $$Q:=\left\{ x\in D(A^r)\::\:\left\Vert A^{r}x\right\Vert _{L_{2}(\mathbb{R})}\leq1\right\} $$
by linear bounded functionals ($A$ is a self-adjoint operator in a Hilbert space $H,$ $f\in H,$ $k<r,\: k,\, r\in\mathbb{N}$), as well as sharp additive inequality of type (\ref{eq:2}), estimating
$\left|F_{f}(x)\right|$ with the help of $\left\Vert x\right\Vert $ and $\left\Vert A^{r}x\right\Vert.$
Even though functional $x^{(k)}(0)$ is not a functional of the described type, Taikov's results follow from results of paper \cite{Babenko_Bilichenko}.

In this paper, we extend the results of \cite{Babenko_Bilichenko} to the case of rather general functions of a self-adjoint operator in a Hilbert space, and present a series of applications.

The paper is organized as follows. In Section  2 we present necessary facts about spectral decomposition of self-adjoint operators and functions of such operators. Section \S 3 contains main results. First, in Section
3.1 we present several technical auxiliary lemmas that are followed by the main results in Section  3.2. We solve the problem of approximation of functionals $(\varphi(A)x, f)$ (where $\varphi(A)$ is some function of self-adjoint operator $A$) on the class of elements of a Hilbert space that is defined with the help of another function $\psi (A)$ of the operator $A$. In addition, in this section we obtain a series of sharp additive inequalities that estimate $|(\varphi(A)x, f)|$ with the help of  $\| \psi (A)x\|$ and $\| x\|$. In Section 3.3  we study properties of functions $N_{\varphi,\psi}$ and $M_{\varphi,\psi}$, which are used in main results. Section 4 is dedicated to applications. In Section  4.1 we find sharp constants in inequalities of the type used in H${\rm{\ddot{o}}}$rmander theorem on comparison of operators (see \cite[Ch. 2, \S 6]{Iosida}) in the case when operators are acting in a Hilbert space and are functions of a self-adjoint operator. In Section 4.2 we generalize Taikov's results to the case of rather arbitrary functions of the differential operator $i\frac d{dt}$. Similar results for operators in the space $L_2({\TT})$ are presented in Section 4.3. Finally, in Section 4.4 we obtain Taikov-type inequalities for integrals with respect to spectral measures, defined with the help of classical orthogonal polynomials.

\section{Spectral decomposition of self-adjoint operators.\\ Functions of operators.}

Following (~\cite[Ch. XIII, \S 1]{Berezanskiy}), we say that there is a partition of unity $E$ defined on a $\sigma$-algebra ${\cal B}$ of Borel sets on $\RR$, if for every $\beta\in {\cal B}$ there is a projector
 $E(\beta)$ in the given Hilbert space $H$ and
\begin{enumerate}
  \item $E(\emptyset)=0,\;\; E(\mathbb{R})=I$;
  \item for any sequence $\{\beta_j\}_{j=1}^\infty\subset {\cal B}$, which consist of mutually disjoint
  sets,
      \[
      E\left(\bigcup\limits_{j=1}^\infty\beta_j\right)=\sum\limits_{j=1}^\infty E(\beta_j).
      \]
\end{enumerate}
With the help of the given partition of unity, one can define so-called spectral integrals (see \cite[Ch. XIII, \S 2]{Berezanskiy}, i.e. integrals of the form
$$
{\cal J}(\varphi)=\displaystyle \int\limits_{\RR} \varphi(t)dE(t),
$$
where  $\varphi \in L_0(\RR)$.
Some necessary properties of such integrals are presented below.
\begin{enumerate}
\item {For any function $\varphi\in L_0(\RR)$} the integral ${\cal J}(\varphi)$ exists as, generally speaking, unbounded operator with dense in $H$ domain
$$
D({\cal J}(\varphi))=\left\{ x\in H\; :\; \displaystyle \int\limits_{\RR} |\varphi(t)|^2d(E(t)x,x)<\infty\right\}.
$$
{If $\varphi\in L_\infty (\RR)$, then this operator is defined on the whole space $H$ and is bounded. Note that for any element $x\in D({\cal J}(\varphi))$
\[
\| {\cal J}(\varphi)x\|^2=\displaystyle \int\limits_{\RR} |\varphi(t)|^2d(E(t)x,x).
\]
}
{If $\varphi\in L_\infty (\RR)$, then for any $x,f\in H$
\[
\left(\displaystyle \int\limits_{\RR}\varphi\left(t\right)dE(t)x,\, f\right)=\left(x,\, \displaystyle \int\limits_{\RR}\overline{\varphi\left(t\right)}dE(t)f\right),
\]
i.e. $({\cal J}(\varphi))^*={\cal J}(\overline{\varphi})$.}
\item For any measurable and almost everywhere finite functions $\varphi$ and $\psi$
\[
\displaystyle \int\limits _{\mathbb{R}}(\varphi(t)+\psi(t))dE(t)=\left(\displaystyle \int\limits _{\mathbb{R}}\varphi(t)dE(t)+\displaystyle \int\limits _{\mathbb{R}}\psi(t)dE(t)\right)^{\sim},
\]
i. e.
\[
{\cal J}(\varphi+\psi)=({\cal J}(\varphi)+{\cal J}(\psi))^\sim ,
\]
and
\[
\displaystyle \int\limits _{\mathbb{R}}\varphi(t)\psi(t)dE(t)=\left(\displaystyle \int\limits _{\mathbb{R}}\varphi(t)dE(t)\displaystyle \int\limits _{\mathbb{R}}\psi(t)dE(t)\right)^{\sim},
\]
i. e.
\[
{\cal J}(\varphi\psi)=({\cal J}(\varphi){\cal J}(\psi))^\sim
\]
{(here, by $A^\sim$ we denote the closure of the operator $A$). In addition, if one of the functions is bounded, then taking closure on the right-hand side of above expressions is no longer needed.}

\end{enumerate}

According to the spectral theorem (see, for instance, ~\cite{Berezanskiy}), for each self-adjoint operator there exists a partition of unity $E$ such that
\begin{equation}\label{sp_th}
A=\displaystyle \int\limits_{\RR}t dE(t),\qquad D(A)=\left\{ x\in H\; :\: \displaystyle \int\limits_{\RR}t^{2}d\left(E(t)x,x\right)<\infty\right\},
\end{equation}
and if $x\in D\left(A\right),$ then
$$
Ax=\displaystyle \int\limits_{\RR}t dE(t)x\qquad \mbox{and}\qquad \left\Vert Ax\right\Vert ^{2}=\displaystyle \int\limits_{\RR}t^{2}d\left(E(t)x,x\right)<\infty.
$$
On the other hand, using (\ref{sp_th}), each partition of unity generates some self-adjoint operator $A$.

Any measurable and almost everywhere finite function $\varphi:\mathbb{R}\to \mathbb{C}$ defines a function $\varphi(A)$ of self-adjoint operator $A$:
\[
\varphi(A)x=\displaystyle \int\limits_{\RR} \varphi(t)dE(t)x={\cal J}_Fx.
\]
Moreover,
\[
D(\varphi(A)):=\left\{ x\in H\; :\; \displaystyle \int\limits_{\RR} |\varphi(t)|^2d(E(t)x,x)<\infty\right\}
\]
and
\[
\| \varphi(A)x\|^2=\displaystyle \int\limits_{\RR} |\varphi(t)|^2d(E(t)x,x).
\]

\medskip

 \section{Main results}

\subsection{Auxiliary lemmas}

Here we present two simple lemmas, which in particular imply that if condition (\ref{condition}) is satisfied, then spectral integrals used further in this paper, exist.

Let functions $\varphi,\psi\in L_0(\RR)$ be given. From now on we assume that these functions satisfy the following condition
\[
\frac{|\varphi|}{1+|\psi|}\in L_\infty(\RR),
\]
or, equivalently,
\begin{equation}\label{condition}
\frac{|\varphi|}{(1+|\psi|^2)^{1/2}}\in L_\infty(\RR).
\end{equation}

\begin{lemma}\label{L1}
If function $\varphi$ and $\psi$ satisfy (\ref{condition}), then
\begin{equation}\label{inclusion}
D(\psi(A))\subset D(\varphi(A)).
\end{equation}
\end{lemma}

{\it Proof.}
 Indeed, let $x\in D(\psi(A))$. Then
\[
\displaystyle \int\limits_{\RR}|\psi(t)|^2d(E(t)x,x)<\infty .
\]
We have
\[
\displaystyle \int\limits_{\RR}|\varphi(t)|^2d(E(t)x,x)=\displaystyle \int\limits_{\RR}\frac{|\varphi(t)|^2}{1+|\psi(t)|^2}(1+|\psi(t)|^2)d(E(t)x,x)
\]
\[
\le \left\|\frac{|\varphi(t)|^2}{1+|\psi(t)|^2}\right\|_{L\infty(\RR)}\displaystyle \int\limits_{\RR}(1+|\psi(t)|^2)d(E(t)x,x)<\infty .
\]
Inclusion (\ref{inclusion}) is proved. $\Box$

\begin{lemma}\label{L2}
If condition (\ref{condition}) is satisfied, then for any $\tau>0$ each of the functions
\[
\frac{\varphi(t)}{1+\tau|\psi(t)|^2},\;\;\frac{\overline{\varphi(t)}}{1+\tau|\psi(t)|^2},\;\;\frac{|\varphi(t)|^2}{(1+\tau|\psi(t)|^2)^2},\;\;
\frac{|\varphi(t)\psi(t)|^2}{(1+\tau|\psi(t)|^2)^2}.
\]
belongs to the space $L_\infty(\RR)$.
\end{lemma}

\medskip

{\it Proof.} We prove the statement for the function $\frac{\varphi(t)}{1+\tau|\psi(t)|^2}$, and the rest can be proved similarly. Without loss of generality, we assume $\tau =1$. The desired statement follows from the inequality
\[
\frac{|\varphi(t)|}{1+|\psi(t)|^2}\le \frac{2|\varphi(t)|}{1+|\psi(t)|},
\]
that is equivalent to
\[
1+|\psi(t)|\le 2+2|\psi(t)|^2\qquad {\rm or}\qquad |\psi(t)|\le 1+2|\psi(t)|^2.
\]
The last inequality for $t$, such that $|\psi(t)|\le 1$, is obvious. Also, for $t$, such that $|\psi(t)|\ge 1$, it follows from the fact that for such $t$ we have $|\psi(t)|\le |\psi(t)|^2$. $\Box$

\subsection{Approximation of functionals\\ and inequalities of Landau-Kolmogorov type}

Let $\varphi,\psi\in L_0(\RR)$, $A$ be a self-adjoint operator in $H$, and
\[
Q_{\psi}:=\left\{ x\in D(\psi(A))\,:\:\left\Vert \psi(A)x\right\Vert \leq1\right\} .
\]
For the element $f\in H$ we define functional $F_{\varphi,f}$ as follows
\[
F_{\varphi,f}(x):=(\varphi(A)x,\, f),\: x\in D(\psi(A)).
\]
In addition, we define functions $N_{\varphi,\psi}(f;\tau)$ and $N_{\varphi,\psi}(f;\tau)$ of $\tau>0$
\[
N_{\varphi,\psi}(f;\tau):=\left\{ \displaystyle \int\limits_{\RR}\frac{\left|\varphi\left(t\right)\right|^{2}}{\left(1+\tau\left|\psi\left(t\right)\right|^{2}\right)^{2}}d\left(E(t)f,f\right)\right\} ^{1/2},
\]
\[
M_{\varphi,\psi}(f;\tau):=\left\{ \displaystyle \int\limits_{\RR}\frac{\left|\varphi(t)\psi\left(t\right)\right|^{2}}{\left(1+\tau\left|\psi\left(t\right)\right|^{2}\right)^{2}}d\left(E(t)f,f\right)\right\} ^{1/2}.
\]

We prove the following two theorems that generalize results of {[}4{]}.
\begin{thm} \label{thm1}
Let functions $\varphi$ and $\psi$ satisfy (\ref{condition}). Then for any $\tau>0$ we have
\begin{equation}
E(N_{\varphi,\psi}(f;\tau))=E(F_{\varphi,f},Q_\psi,N_{\varphi,\psi}(f;\tau))=\tau M{}_{\varphi,\psi}(f;\tau).\label{eq:15}
\end{equation}
In addition, the extremal approximating functional for $F_{\varphi,f}$ is the functional
\begin{equation}\label{extr_funct}
g_{\tau}(x):=\displaystyle \int\limits_{\RR}\frac{\varphi\left(t\right)}{1+\tau\left|\psi\left(t\right)\right|^{2}}d\left(E(t)x,f\right).
\end{equation}
\end{thm}

\begin{thm}\label{thm2}
Under assumptions of Theorem \ref{thm1}, for any $x\in D(\psi(A))$ and any
$\tau>0$ we have
\begin{equation}
|F_{\varphi,f}(x)|\leq\tau M{}_{\varphi,\psi}(f;\tau)\left\Vert \psi(A)x\right\Vert +N_{\varphi,\psi}(f;\tau)\left\Vert x\right\Vert .\label{eq:16}
\end{equation}
Inequality (\ref{eq:16}) becomes equality for the element
\begin{equation}\label{xt}
x_{\tau}:=\displaystyle \int\limits_{\RR}\frac{\overline{\varphi\left(t\right)}}{1+\tau\left|\psi\left(\tau\right)\right|^{2}}dE(t)f.
\end{equation}
\end{thm}

{\it Proof.} Taking into account the general scheme of solving Stechkin's problem that was presented in Introduction, proofs of Theorems \ref{thm1} and \ref{thm2} consist of the following four steps.
\begin{enumerate}
\item First we prove that
\begin{equation}
E(N_{\varphi,\psi}(f;\tau))\le\tau M{}_{\varphi,\psi}(f;\tau).\label{up_est}
\end{equation}
\item Estimate (\ref{up_est}) implies inequality (\ref{eq:16}).
\item Next, we prove that inequality (\ref{up_est}) becomes equality for the element defined in (\ref{xt}).
\item Using this fact, we prove equality (\ref{eq:15}). Having that, both theorems are proved.
\end{enumerate}

{\bf Step 1.}
First let us show that for functionals $g_\tau$ (see (\ref{extr_funct}))
\begin{equation}\label{norm_est}
\left\Vert g_{\tau}\right\Vert \leq N_{\varphi,\psi}(f;\tau).
\end{equation}
 Indeed,
\begin{align*}
\left|g_{\tau}(x)\right| &=
\left|\displaystyle \displaystyle \int\limits_{\RR}\frac{\varphi\left(t\right)}{1+\tau\left|\psi\left(t\right)\right|^{2}}d\left(E(t)x,f\right)\right|=\left|\left(\displaystyle \displaystyle \int\limits_{\RR}\frac{\varphi\left(t\right)}{1+\tau\left|\psi\left(t\right)\right|^{2}}dE(t)x,f\right)\right|\\
&=
\left|\left(x,\displaystyle  \displaystyle \int\limits_{\RR}\frac{\overline{\varphi\left(t\right)}}{1+\tau\left|\psi\left(t\right)\right|^{2}}dE(t)f\right)\right|
\le   \| x\| \cdot \left\|\displaystyle  \displaystyle \int\limits_{\RR}\frac{\overline{\varphi\left(t\right)}}{1+\tau\left|\psi\left(t\right)\right|^{2}}dE(t)f\right\| \\
&= \| x\| \cdot \left\{\displaystyle \displaystyle \int\limits_{\RR}\frac{\left|\varphi\left(t\right)\right|^{2}}{\left(1+\tau\left|\psi\left(t\right)\right|^{2}\right)^{2}}d\left(E(t)f,f\right)\right\} ^{1/2}=N_{\varphi,\psi}(f;\tau)\left\Vert x\right\Vert.
\end{align*}
This gives the desired estimate (\ref{norm_est}).

If  $x\in Q_{\psi},$  we obtain
\begin{align*}
\left|F_{\varphi,f}(x)-g_{\tau}(x)\right|&=\left| \displaystyle \displaystyle \int\limits_{\RR}\varphi\left(t\right)d\left(E(t)x,f\right)-\displaystyle \int\limits_{\RR}\frac{\varphi\left(t\right)}{1+\tau\left|\psi\left(t\right)\right|^{2}}d\left(E(t)x,f\right)\right| \\
&=\tau\left|\displaystyle \int\limits_{\RR}\frac{\varphi\left(t\right)\overline{\psi(t)}}{1+\tau\left|\psi\left(t\right)\right|^{2}}{\psi(t)}d\left(E(t)x,f\right)\right| \\
&=\tau\left|\left(\displaystyle \int\limits_{\RR}\frac{\varphi\left(t\right)\overline{\psi(t)}}{1+\tau\left|\psi\left(t\right)\right|^{2}}{\psi(t)}dE(t)x,f\right)\right| \\
&=\tau\left|\left(\displaystyle \int\limits_{\RR}{\psi(t)}dE(t)x,\displaystyle \int\limits_{\RR}\frac{\overline{\varphi\left(t\right)}{\psi(t)}}{1+\tau\left|\psi\left(t\right)\right|^{2}}dE(t)f\right)\right| \\
&\le  \tau\left\|\displaystyle \int\limits_{\RR}{\psi(t)}dE(t)x\right\|\cdot\left\|\displaystyle \int\limits_{\RR}\frac{\overline{\varphi\left(t\right)}{\psi(t)}}{1+\tau\left|\psi\left(t\right)\right|^{2}}dE(t)f\right\| \\
&=\tau\left\|{\psi}(A)x\right\|\cdot\left\|\displaystyle \int\limits_{\RR}\frac{\overline{\varphi\left(t\right)}{\psi(t)}}{1+\tau\left|\psi\left(t\right)\right|^{2}}dE(t)f\right\|\\
&=\tau \left\|{\psi}(A)x\right\| \cdot \left(\displaystyle \int\limits_{\RR}\frac{\left|\varphi\left(t\right)\psi(t)\right|^{2}}{\left(1+\tau\left|\psi\left(t\right)\right|^{2}\right)^{2}}d\left(E(t)f,f\right)\right)^{1/2}\le\tau M{}_{\varphi,\psi}(f;\tau).
\end{align*}
Thus, for any $x\in Q_\psi$
\begin{equation}
\left|F_{\varphi,f}(x)-g_{\tau}(x)\right|\leq\tau M{}_{\varphi,\psi}(f;\tau)\label{eq:18}
\end{equation}
and, hence,
\begin{equation}
E(N_{\varphi,\psi}(f;\tau))\leq\tau M{}_{\varphi,\psi}(f;\tau)\label{eq:19}
\end{equation}
Inequality (\ref{up_est}) is proved.

{\bf Step 2.}
From estimates (\ref{eq:18}) and (\ref{eq:19}), it follows that for
$\tau>0$ and any $x\in D(\psi(A))$
\begin{align*}
\left|F_{\varphi,f}(x)\right| & \leq  \left|F_{\varphi,f}(x)-g_{\tau}(x)\right|+\left|g_{\tau}(x)\right| \\
&=\left|F_{\varphi,f}\left(\frac x{\|\psi(A)x\|}\right)-g_{\tau}\left(\frac x{\|\psi(A)x\|}\right)\right|\|\psi(A)x\|+\left\|g_{\tau}\right\|\| x\| \\
&\leq \tau M{}_{\varphi,\psi}(f;\tau)\left\Vert \psi(A)x\right\Vert +N_{\varphi,\psi}(f;\tau)\left\Vert x\right\Vert.
\end{align*}
Inequality (\ref{eq:16}) is proved.

{\bf Step 3.}
By Lemma \ref{L2}, $\frac{\overline{\varphi}\psi}{1+\tau |\psi|^2}\in L_\infty(\RR)$, which implies that the element $x_{\tau},$ defined by equality (\ref{xt}), belongs to $x_\tau\in D(\psi(A))\subset D(\varphi(A))$. Let us show that for this element inequality (\ref{eq:16}) becomes an equality. We have
\begin{equation}
\left\Vert x_{\tau}\right\Vert =\left(\displaystyle \int\limits_{\RR}\frac{\left|\varphi\left(t\right)\right|^{2}}{\left(1+\tau\left|\psi\left(t\right)\right|^{2}\right)^{2}}d\left(E(t)f,f\right)\right)^{1/2};\label{eq:22}
\end{equation}
\[
\varphi(A)x_{\tau}=\displaystyle \int\limits_{\RR}\frac{\left|\varphi\left(t\right)\right|^{2}}{1+\tau\left|\psi\left(t\right)\right|^{2}}dE(t)f,
\]
\[
\psi(A)x_{\tau}=\displaystyle \int\limits_{\RR}\frac{\overline{\varphi\left(t\right)}\psi(t)}{1+\tau\left|\psi\left(t\right)\right|^{2}}dE(t)f;
\]
\begin{equation}
\left\Vert \psi(A)x_{\tau}\right\Vert =\left(\displaystyle \int\limits_{\RR}\frac{\left|\varphi\left(t\right)\psi(t)\right|^{2}}{\left(1+\tau\left|\psi\left(t\right)\right|^{2}\right)^{2}}d\left(E(t)f,f\right)\right)^{1/2};\label{eq:25}
\end{equation}
\begin{equation}
|F_{\varphi,f}(x_\tau)|=\left|\left(\varphi(A)x_{\tau},\, f\right)\right|=\displaystyle \int\limits_{\RR}\frac{\left|\varphi\left(t\right)\right|^{2}}{1+\tau\left|\psi\left(t\right)\right|^{2}}d\left(E(t)f,f\right).\label{eq:26}
\end{equation}
Substituting (\ref{eq:26}), (\ref{eq:22}), and (\ref{eq:25}) in (\ref{eq:16}),
we obtain
\begin{align*}
|F_{\varphi,f}(x_\tau)|&=\displaystyle \displaystyle \int\limits_{\RR}\frac{\left|\varphi\left(t\right)\right|^{2}}{1+\tau\left|\psi\left(t\right)\right|^{2}}d\left(E(t)f,f\right) \\
&=\displaystyle \int\limits_{\RR}\frac{\left|\varphi\left(t\right)\right|^{2}\left(1+\tau\left|\psi\left(t\right)\right|^{2}\right)}{\left(1+\tau\left|\psi\left(t\right)\right|^{2}\right)^{2}}d\left(E(t)f,f\right) \\
&=\displaystyle \int\limits_{\RR}\frac{\left|\varphi\left(t\right)\right|^{2}}{\left(1+\tau\left|\psi\left(t\right)\right|^{2}\right)^{2}}d\left(E(t)f,f\right)+\tau\displaystyle \int\limits_{\RR}\frac{\left|\varphi\left(t\right)\psi(t)\right|^{2}}{\left(1+\tau\left|\psi\left(t\right)\right|^{2}\right)^{2}}d\left(E(t)f,f\right)\\
&=N_{\varphi,\psi}(f;\tau)\left\Vert x_{\tau}\right\Vert +\tau M{}_{\varphi,\psi}(f;\tau)\left\Vert \psi(A)x_{\tau}\right\Vert .
\end{align*}
Thus,
\[
|F_{\varphi,f}(x_\tau)|=N_{\varphi,\psi}(f;\tau)\left\Vert x_{\tau}\right\Vert +\tau M{}_{\varphi,\psi}(f;\tau)\left\Vert \psi(A)x_{\tau}\right\Vert ,
\]
and, hence, for $y_{\tau}=\frac{x_{\tau}}{\left\Vert \psi(A)x_{\tau}\right\Vert }$ we have
\begin{equation}\label{equal}
|F_{\varphi,f}(y_\tau)|=\tau M{}_{\varphi,\psi}(f;\tau)+N_{\varphi,\psi}(f;\tau)\left\Vert y_{\tau}\right\Vert .
\end{equation}

{\bf Step 4.} We now prove
\[
E(N_{\varphi,\psi}(f;\tau))=\tau M{}_{\varphi,\psi}(f;\tau).
\]

For the element $y_\tau$ we have equality (\ref{equal}).
Besides that,
\[
\left| F_{\varphi,f}(y_\tau)\right|\leq E\left(N_{\varphi,\psi}(f;\tau)\right)+N_{\varphi,\psi}(f;\tau)\left\Vert y_{\tau}\right\Vert .
\]
Comparing the last inequality and (\ref{equal}), we see that
\[
E\left(N_{\varphi,\psi}(f;\tau)\right)\geq\tau M{}_{\varphi,\psi}(f;\tau).
\]
Together with (\ref{eq:19}), it gives
\[
E\left(N_{\varphi,\psi}(f;\tau)\right)=\tau M{}_{\varphi,\psi}(f;\tau).
\]
Theorems \ref{thm1} and \ref{thm2} are now proved. $\Box$

\subsection{Properties of functions $N_{\varphi,\psi}(f,\tau)$ and $M_{\varphi,\psi}(f,\tau)$.}

Because of Theorems \ref{thm1} and \ref{thm2}, it is interesting to study further the question under what conditions on functions $\varphi$ and $\psi$ one could claim that the problem of computing the quantity $E\left(F_{\varphi,\psi},Q_\psi,N\right)$ is solved for all $N\in (0,\|\varphi(A)\|)$ in the case $f\in D(\varphi(A))$, and for all $N\in (0,\infty)$ in the case $f\notin D(\varphi(A))$.
Lemmas below provide sufficient conditions.

\begin{lemma}\label{L3}
We assume functions $\varphi$ and $\psi$ satisfy condition (\ref{condition}). Then function $N_{\varphi,\psi}(f;\tau)$
continuously depends on $\tau$ and is non-increasing with $\tau$.
\end{lemma}

{\it Proof.} The fact that $N_{\varphi,\psi}(f;\tau)$ is non-increasing is obvious. Continuity at every point $\tau \in (0,+\infty)$ readily follows from the following equality that holds for all $0<\tau_1<\tau_2$
\[
N_{\varphi,\psi}(f;\tau_1)^2-N_{\varphi,\psi}(f;\tau_2)^2
\]
\[
=(\tau_2-\tau_1)\displaystyle \int\limits_{\RR}
\frac{|\varphi(t)\psi(t)|^2(2+(\tau_1+\tau_2)|\psi(t)|^2)}{(1+\tau_1|\psi(t)|^2)^2(1+\tau_2|\psi(t)|^2)^2}d(E(t)f,f)
\]
and from the fact that functions $\frac{|\varphi(t)\psi(t)|^2}{(1+\tau_1|\psi(t)|^2)^2}$ and $\frac{2+(\tau_1+\tau_2|\psi(t)|^2}{1+\tau_2|\psi(t)|^2)^2}$ are bounded. $\Box$

Next we study the behavior of this function as $\tau\to 0$.

\begin{lemma}\label{L4}
Let functions $\varphi$ and $\psi$ satisfy condition (\ref{condition}) and be continuous. Then
$$
\lim\limits_{\tau\to 0}N_{\varphi,\psi}(f;\tau)=\left\{\begin{array}{clcr}
\mbox{$\|\varphi(A)f\|$}, & \mbox{ if $f\in D(\varphi(A)),$} \\
\mbox{$+\infty$} , & \mbox{ if $f\notin D(\varphi(A))$}.
\end{array}\right.
$$
\end{lemma}

{\it Proof.} Let $f\in D(\varphi(A))$. Then
 \[
 \displaystyle \int\limits_{\RR}|\varphi(t)|^2d(E(t)f,f)=\|\varphi(A)f\|^2.
 \]
 For any $\varepsilon>0$ there exists $a>0$ such that
 \[
 \displaystyle \int_{-a}^a|\varphi(t)|^2d(E(t)f,f)>(1-\varepsilon)\|\varphi(A)f\|^2.
 \]
Since functions $\varphi$ and $\psi$ are continuous, for all small enough $\tau$ we have
\begin{align*}
\|\varphi(A)f\|^2&\ge \displaystyle \int\limits_{\RR}\frac{|\varphi(t)|^2}{(1+\tau|\psi(t)|^2)^2}d(E(t)f,f)
\ge\displaystyle \int_{-a}^a\frac{|\varphi(t)|^2}{(1+\tau|\psi(t)|^2)^2}d(E(t)f,f)\\
&> \displaystyle \int_{-a}^a(1-\varepsilon)|\varphi(t)|^2d(E(t)f,f)
>(1-\varepsilon)^2\|\varphi(A)f\|^2.
\end{align*}
In the considered case, the desired fact is proved.

 Let now $f\notin D(\varphi(A))$. Then for any $N>0$ there exists $a>0$ such that
 \[
 \displaystyle \int_{-a}^a|\varphi(t)|^2d(E(t)f,f)>N
 \]
and for all small enough $\tau$ we have
\[
\displaystyle \int\limits_{\RR}\frac{|\varphi(t)|^2}{(1+\tau|\psi(t)|^2)^2}d(E(t)f,f)
\ge\displaystyle \int_{-a}^a\frac{|\varphi(t)|^2}{(1+\tau|\psi(t)|^2)^2}d(E(t)f,f)>\frac N2.
\]
Lemma is proved. $\Box$

\begin{lemma}\label{L5}
Let functions $\varphi$ and $\psi$ satisfy (\ref{condition}), be continuous, and
\begin{equation}\label{condition3}
\sup\limits_{t\in\RR}\frac{|\varphi(t)|^2}{(1+\tau|\psi(t)|^2)^2}\to 0,\qquad {\rm as}\qquad \tau\to\infty .
\end{equation}
Then
\begin{equation}\label{limitus}
\lim\limits_{\tau\to \infty}N_{\varphi,\psi}(f;\tau)=0,
\end{equation}
and, hence, for any $N\in (0,\|\varphi(A)f\|)$ (for any $N\in (0,+\infty)$ in the case $f\notin D(\varphi(A))$) equation for $\tau$
\[
N_{\varphi,\psi}(f;\tau)=N
\]
has a solution.
\end{lemma}

{\it Proof} By Lemmas \ref{L3} and \ref{L4}, together with properties of continuous functions, the statement of the lemma is proved once we prove (\ref{limitus}), which is obvious if condition (\ref{condition3}) is satisfied. $\Box$

{\bf Remark 1.}
In order to condition (\ref{condition3}) hold, it is sufficient to have
\begin{equation}\label{condition4}
\sup\limits_{t\in\RR}\frac{|\varphi(t)|^2}{1+\tau|\psi(t)|^2}\to 0,\qquad {\rm as}\qquad \tau\to\infty .
\end{equation}

{\bf Remark 2.}
It is easy to verify that condition (\ref{condition4}) is satisfied for functions $\varphi(t)=t^\alpha$ and $\psi(t)=t^\beta$, $0<\alpha<\beta$.

\medskip

\begin{lemma}\label{L6}
Let functions $\varphi$ and $\psi$ satisfy condition (\ref{condition}) and be continuous. In addition, let either $f\in D(\varphi(A))$, or {$f\notin D(\varphi(A))$ and}
\begin{equation}\label{condition5}
\tau\sup\limits_{t\in\RR}\frac{|\varphi(t)|^2}{1+\tau |\psi(t)|^2}\to 0, \qquad {\rm as}\qquad \tau \to 0.
\end{equation}
Then
\begin{equation}\label{stat}
\tau M_{\varphi,\psi}(f;\tau)\to 0, \qquad {\rm as}\qquad \tau \to 0.
\end{equation}
\end{lemma}

{\it Proof.} When $f\in D(\varphi(A))$ the function $M_{\varphi,\psi}(f;\tau)$ is bounded and (\ref{stat}) is obvious. In the case when $f\notin D(\varphi(A))$, the statement follows from the obvious estimate
\[
M_{\varphi,\psi}(f;\tau)^2\le \max\limits_{t\in\RR}\frac{|\varphi(t)|^2}{1+\tau |\psi(t)|^2}\| f\|^2
\]
and (\ref{condition5}). $\Box$

{\bf Remark 3.}
Note that condition (\ref{condition5}) is satisfied, for instance, by functions $\varphi(t)=t^\alpha$ and $\psi(t)=t^\beta$, $0<\alpha<\beta$, from Remark 2.

\section{Applications.}

\subsection{Sharp constant in H${\bf{\rm{\ddot{o}}}}$rmander theorem\\ on comparison of operators in the case of Hilbert spaces.}

The next theorem is due to H${\rm{\ddot{o}}}$rmander (see, for instance, \cite[Ch. 2, \S  6, p. 117]{Iosida}).
\begin{thm}\label{thm3}
Let us consider Banach spaces $X_i\; (i=0,1,2;\; X_0=X)$ and linear operators $T_i \; (i=1,2)$, that map $D(T_i)\subset X$ into spaces $X_i$. Let operator $T_1$ be closed and let operator $T_2$ admit closed extension. If $D(T_1)\subset D(T_2)$, then there exists a constant $C$, such that for all $x\in D(T_1)$
\begin{equation}\label{herm1}
\| T_2x\|_{X_2}\le C\{ \| x\|_{X}^2+\| T_1x\|_{X_1}^2\}^{1/2}.
\end{equation}
\end{thm}

From this theorem it follows that for any functional $f\in X_2^*$ and any $\tau>0$ there exists a constant $C(f,\tau)$ such that for all $x\in D(T_1)$
\begin{equation}\label{herm2}
( T_2x,f)\le C(f,\tau)\{ \| x\|_{X}^2+\tau\| T_1x\|_{X_1}^2\}^{1/2}.
\end{equation}

In the case when $X_0=X_1=X_2=H$ and operators $T_1$ and $T_2$ are functions of a self-adjoint operator in $H$, we find sharp constants in inequalities (\ref{herm2}) and (\ref{herm1}).

\begin{thm}\label{thm4}
Let functions $\phi$ and $\psi$ satisfy (\ref{condition}), and let $A$ be a self-adjoint operator in a Hilbert space $H$. Then for any $f\in H$, arbitrary $\tau>0$, and arbitrary $x\in D(\psi(A))$ we have
\begin{equation}\label{herm}
\left|\left(\varphi(A)x,\, f\right)\right|\leq \left\{\displaystyle \int\limits_{\RR}\frac{\left|\varphi\left(t\right)\right|^{2}}{1+\tau\left|\psi\left(t\right)\right|^{2}}d\left(E(t)f,f\right)
\right\}^{1/2}\{\left\Vert x\right\Vert^2 +\tau\left\Vert \psi(A)x\right\Vert^2\}^{1/2} .
\end{equation}
Inequality (\ref{herm}) becomes an equality for the element defined in (\ref{xt}).
\end{thm}

This is precisely the sharp form of the corollary from H${\rm{\ddot{o}}}$rmander's theorem that deals with functionals. Note that this inequality  (\ref{herm}) for powers of a self-adjoint operator $A$  was proved in \cite{babenko_bilichenko}.

{\it Proof.} Applying Schwartz  inequality to the right-hand side of (\ref{eq:16}), we obtain
\[
\left|\left(\varphi(A)x,\, f\right)\right|\leq \{N_{\varphi,\psi}(f;\tau)^2+\tau M{}_{\varphi,\psi}(f;\tau)^2\}^{1/2}\{\left\Vert x\right\Vert^2 +\tau\left\Vert \psi(A)x\right\Vert^2\}^{1/2} .
\]
Since
\[
N_{\varphi,\psi}(f;\tau)^2+\tau M{}_{\varphi,\psi}(f;\tau)^{2}
\]
\[
=\displaystyle \int\limits_{\RR}\frac{\left|\varphi\left(t\right)\right|^{2}}{\left(1+\tau\left|\psi\left(t\right)\right|^{2}\right)^{2}}d\left(E(t)f,f\right)+\tau\displaystyle \int\limits_{\RR}\frac{\left|\varphi\left(t\right)\psi(t)\right|^{2}}{\left(1+\tau\left|\psi\left(t\right)\right|^{2}\right)^{2}}d\left(E(t)f,f\right)
\]
\[
=\displaystyle \int\limits_{\RR}\frac{\left|\varphi\left(t\right)\right|^{2}}{1+\tau\left|\psi\left(t\right)\right|^{2}}d\left(E(t)f,f\right),
\]
we arrive at inequality (\ref{herm}). The fact that inequality (\ref{herm}) becomes an equality for $x_\tau$, defined by (\ref{xt}), can be verified similar to the way it was done in Theorem \ref{thm2}. $\Box$

Next we obtain H${\rm{\ddot{o}}}$rmander's theorem in a Hilbert space. The obtained inequalities can naturally be called Hardy-Littlewood-Polya-type inequalities. More information about other Hardy-Littlewood-Polya-type inequalities can be found in, for instance, works \cite{babenko_inequalities},  \cite{Babenko_Babenko_Kriachko}, which also provide further references.

{
\begin{thm}\label{HLP} Let continuous on $\RR$ functions $\phi$ and $\psi$ satisfy (\ref{condition}), and $A$ be a self-adjoint operator in a Hilbert space $H$. Then for any $f\in H$, arbitrary $\tau>0$ and arbitrary $x\in D(\psi(A))$, we have
\begin{equation}\label{104}
\|\varphi(A)x\|\le\sup\limits_{t\in \mathbb{R}}\left\{\frac{| \varphi(t)|^2}{1+\tau| \psi(t)|^2}\right\}^{\frac 12}\left\{ \left\Vert x\right\Vert ^2+\tau\left\Vert \psi(A)x\right\Vert ^{2}\right\} ^{\frac{1}{2}}.
\end{equation}
If partition of unity $E(\beta)$ is such that for any $t\in \mathbb{R}$ and $\delta>0$, we have
\begin{equation}\label{condition2}
E([t,t+\delta])\neq 0,
\end{equation}
then inequality (\ref{104}) is sharp.
\end{thm}
}
{\it Proof.}
 From inequality (\ref{herm}) we obtain
\[
\|(\varphi\left(A\right)x\|=\sup\limits_{f\in H\atop \| f\|\le 1}\left|\left(\varphi\left(A\right)x,\, f\right)\right|
\]
\[
\leq \sup\limits_{f\in H\atop \| f\|\le 1}\left\{\displaystyle \int\limits_{\RR}\frac{\left|\varphi\left(t\right)\right|^{2}}{1+\tau\left|\psi\left(t\right)\right|^{2}}d\left(E(t)f,f\right)
\right\}^{1/2}\{ \left\Vert x\right\Vert ^2+\tau\left\Vert \psi\left(A\right)x\right\Vert ^{2}\} ^{1/2}.
\]
We have
\[
\sup\limits_{f\in H\atop \| f\|\le 1}\displaystyle \int\limits_{\RR}\frac{\left|\varphi\left(t\right)\right|^{2}}{1+\tau\left|\psi\left(t\right)\right|^{2}}d\left(E(t)f,f\right)
\]
\[
\le \sup\limits_{f\in H\atop \| f\|\le 1}\sup\limits_{t\in \mathbb{R}}\frac{| \varphi (t)|^2}{1+\tau| \psi (t)|^2}\displaystyle \int\limits_{\RR} d(E(t)f,f)=\sup\limits_{t\in \mathbb{R}}\frac{| \varphi (t)|^2}{1+\tau|\psi (t)|^2}.
\]
Inequality (\ref{104}) is proved.

We show that if (\ref{condition2}) is satisfied, then the constant in this  inequality is the best possible.

We assume that there exists $\varepsilon>0$ such that for $x\in D(\psi(A))$
\begin{equation}\label{contra}
\|\varphi(A)x\|^2\le(1-\varepsilon)\sup\limits_{t\in \mathbb{R}}\frac{| \varphi(t)|^2}{1+\tau| \psi(t)|^2}\left\{ \left\Vert x\right\Vert ^2+\tau\left\Vert \psi(A)x\right\Vert ^{2}\right\}.
\end{equation}
Taking into account continuity of $\varphi$ and $\psi$ we find a point $t_\varepsilon$ and a number $\delta>0$ such that
\begin{enumerate}
\item $\sup\limits_{t\in \mathbb{R}}\frac{| \varphi (t)|^2}{1+\tau| \psi (t)|^2}<\sqrt{1+\varepsilon}\frac{| \varphi (t_\varepsilon)|^2}{1+\tau| \psi (t_\varepsilon)|^2}$;
\item $|\varphi(t)|^2>(1-\varepsilon^2)|\varphi(t_\varepsilon)|^2$ for any $t\in [t_\varepsilon, t_\varepsilon+\delta]$;
\item $1+\tau\left|\psi\left(t\right)\right|^{2}<\sqrt{1+\varepsilon}(1+\tau\left|\psi\left(t_\varepsilon\right)\right|^{2})$ for any $t\in [t_\varepsilon, t_\varepsilon+\delta]$.
\end{enumerate}

Further, we choose $f_\varepsilon\in E([t_\varepsilon , t_\varepsilon +\delta])(H)$ so that $\| f_\varepsilon\|=1$. Applying properties 1-3 and inequality (\ref{contra}) for the element $f_\varepsilon$, we obtain
\begin{align*}
(1-\varepsilon^2)|\varphi(t_\varepsilon)|^2&<\displaystyle \int_{t_\varepsilon}^{t_\varepsilon+\delta}|\varphi(t)|^2d(E(t)f_\varepsilon,f_\varepsilon)=
\| \varphi(A)f_\varepsilon\|^2 \\
&\le  (1-\varepsilon)\sup\limits_{t\in \mathbb{R}}\frac{| \varphi(t)|^2}{1+\tau| \psi(t)|^2}\left\{ \left\Vert f_\varepsilon\right\Vert ^2+\tau\left\Vert \psi(A)f_\varepsilon\right\Vert ^{2}\right\} \\
&\le  (1-\varepsilon)\sqrt{1+\varepsilon}\cdot\frac{| \varphi (t_\varepsilon)|^2}{1+\tau| \psi (t_\varepsilon)|^2}\left(1+\tau \displaystyle \int_{t_\varepsilon}^{t_\varepsilon+\delta}|\psi(t)|^2d(E(t)f_\varepsilon,f_\varepsilon)\right) \\
&=(1-\varepsilon)\sqrt{1+\varepsilon}\cdot\frac{| \varphi (t_\varepsilon)|^2}{1+\tau| \psi (t_\varepsilon)|^2} \displaystyle \int_{t_\varepsilon}^{t_\varepsilon+\delta}(1+\tau|\psi(t)|^2)d(E(t)f_\varepsilon,f_\varepsilon)\\
&\le(1-\varepsilon)\sqrt{1+\varepsilon}\cdot\frac{| \varphi (t_\varepsilon)|^2}{1+\tau| \psi (t_\varepsilon)|^2} \sqrt{1+\varepsilon}\cdot(1+\tau|\psi(t_\varepsilon)|^2)=(1-\varepsilon^2)| \varphi (t_\varepsilon)|^2.
\end{align*}
This brings us to a contradictory inequality
\[
(1-\varepsilon^2)<(1-\varepsilon^2).
\]
Theorem is proved. $\Box$

\subsection{Generalization of Taikov's results.}

As usual, by Fourier transform of a function $x\in L_2(\RR)$, we understand
$$
{\cal{F}}[x](t):=\frac{1}{\sqrt{2\pi}}\displaystyle \displaystyle \int\limits_{\RR}e^{-its}x(s)ds=\frac{1}{\sqrt{2\pi}}\lim\limits_{N\to \infty}\displaystyle \displaystyle \int_{-N}^Ne^{-its}x(s)ds,
$$
and by inverse Fourier transform we understand
$$
{\cal{F}}^{-1}[y](t):=\frac{1}{\sqrt{2\pi}}\displaystyle \displaystyle \int\limits_{\RR}e^{its}y(s)ds=\frac{1}{\sqrt{2\pi}}\lim\limits_{N\to \infty}\displaystyle \displaystyle \int_{-N}^Ne^{its}y(s)ds.
$$
Let us consider operator $A=\frac 1i\frac{d}{dt}$ in $L_2(\RR)$. From properties of inverse Fourier transform, it follows that for the partition of unity corresponding to this operator the following is true:
\begin{equation}\label{idt}
E([\alpha,\beta])x={\cal{F}}^{-1}[{\cal{F}}[x] \chi_{[\alpha,\beta]}]=\frac{1}{2\pi}\displaystyle \int\limits_{\RR}\frac{e^{i\beta(z-u)}-e^{i\alpha(z-u)}}{i(z-u)}x(z)dz,\;\alpha<\beta.
\end{equation}
For the function $\varphi(A)$ of the operator $A$ (here and below functions $\varphi,\psi$ are continuous and satisfy (\ref{condition})) we compute the spectral integral and write
\[
\varphi(A)x(t)=\left(\displaystyle \displaystyle \int\limits_{\RR}\varphi(s) dE(s)x\right)(t)={\cal{F}}^{-1}[{\varphi\cal{F}}[x]](t).
\]
The fact that function $x\in L_2(\RR)$ belongs to the domain of the operator $\varphi(A)$ implies $\varphi {\cal{F}}[x]\in L_2(\RR)$. If $x\in D(\psi(A))$ and $\frac{|\varphi|}{(1+|\psi|^2)^{1/2}}\in L_2(\RR)$, then it is easy to see that $\varphi {\cal{F}}[x]\in L_1(\RR)$, and, hence, the function $\varphi(A)x(t)$ is continuous.

Next, for any $a>0$ we set
$$
f_a(t):={\cal{F}}^{-1}[\chi_{[-a,a]}](t).
$$
Using the fact that operator ${\cal{F}}^{-1}$ is unitary, we have
\begin{align*}
|(\varphi(A)x,f_a)|&=\left( {\cal{F}}^{-1}[{\varphi \cal{F}}[x]], f_a\right)=\left( {\varphi \cal{F}}[x],{\cal{F}}[f_a]\right)\\
&=\left( {\varphi \cal{F}}[x],\chi_{[-a,a]}\right)=\displaystyle \displaystyle \int_{-a}^a \varphi(s){ \cal{F}}[x](s)ds.
\end{align*}
As $a\to +\infty$, we have
$$
\displaystyle \displaystyle \int_{-a}^a \varphi(s){ \cal{F}}[x](s)ds \to \displaystyle \displaystyle \int\limits_{\RR} \varphi(s){ \cal{F}}[x](s)ds ={\cal{F}}^{-1}[{\varphi \cal{F}}[x]](0)=\varphi(A)x(0).
$$

From the fact that $\frac{|\varphi|}{(1+|\psi|^2)^{1/2}} \in L_2(\RR)$, it follows that $\frac{|\varphi|}{1+|\psi|^2}$ and $\frac{|\varphi\psi|}{1+|\psi|^2}\in L_2(\RR)$. Then for $a\to \infty$
\begin{align*}
N_{\varphi,\psi}(f_a;\tau)^2&:= \displaystyle \int\limits_{\RR}\frac{\left|\varphi\left(t\right)\right|^{2}}{\left(1+\tau\left|\psi\left(t\right)\right|^{2}\right)^{2}}d\left(E(t)f_a,f_a\right) \\
&=\left ( {\cal{F}}^{-1}\left [ \frac{|\varphi |^2}{(1+\tau |\psi|^2)^2}{\cal{F}}[f_a]\right ], f_a\right )\\
&=\left ( \frac{|\varphi |^2}{(1+\tau |\psi|^2)^2}\chi_{[-a,a]}, \chi_{[-a,a]}  \right )\\
&= \displaystyle \displaystyle \int_{-a}^{a}\frac{|\varphi(s)|^2}{(1+\tau|\psi(s)|^2)^2}ds \to  \displaystyle \displaystyle \int\limits_{\RR}\frac{|\varphi(s)|^2}{(1+\tau|\psi(s)|^2)^2}ds,
\end{align*}
and
\[
M{}_{\varphi,\psi}(f_a;\tau)^2\to \displaystyle \int\limits_{\RR}\frac{\left|\varphi\left(s\right)\psi(s)\right|^{2}}{\left(1+\tau\left|\psi\left(s\right)\right|^{2}\right)^{2}}ds. \]

Thus, taking the limit as $a\to \infty$ in the inequality
$$
|(\varphi(A)x,f_a)|\leq N_{\varphi,\psi}(f_a;\tau) \|x\|_{L_2(\RR)}+\tau M{}_{\varphi,\psi}(f_a;\tau)\|\psi(A)x\|_{L_2(\RR)},
$$
we arrive at
$$
|\varphi(A)x(0)|\leq \left ( \displaystyle \displaystyle \int\limits_{\RR}\frac{|\varphi(s)|^2}{(1+\tau|\psi(s)|^2)^2}ds \right )^{1/2}\|x\|_{L_2(\RR)}
 $$
 \begin{equation}\label{taik_gen_1}
\qquad\qquad\qquad\qquad +\tau \left ( \displaystyle \displaystyle \int\limits_{\RR}\frac{|\varphi(s)\psi(s)|^2}{(1+\tau|\psi(s)|^2)^2}ds \right )^{1/2}\|\psi(A)x\|_{L_2(\RR)}.
 \end{equation}
Using the fact that the operators $\varphi(A)$ and $\psi(A)$ are invariant with respect to the shift, the last inequality implies
$$
\|\varphi(A)x\|_{L_{\infty}(\RR)}\leq \left ( \displaystyle \displaystyle \int\limits_{\RR}\frac{|\varphi(s)|^2}{(1+\tau|\psi(s)|^2)^2}ds \right )^{1/2}\|x\|_{L_2(\RR)}
 $$
  \begin{equation}\label{taik_gen_2}
 \qquad\qquad\qquad\qquad
 +\tau \left ( \displaystyle \displaystyle \int\limits_{\RR}\frac{|\varphi(s)\psi(s)|^2}{(1+\tau|\psi(s)|^2)^2}ds \right )^{1/2}\|\psi(A)x\|_{L_2(\RR)}.
 \end{equation}
From the last inequality one can easily obtain Taikov's inequality by taking $\varphi(t)=t^k,\; \psi(t)=t^r,\;\; k,r\in\NN, \; k<r$.

Therefore, we proved the following theorem
\begin{thm}
Assume that functions $\varphi,\psi$ are continuous, satisfy condition (\ref{condition}), and $\frac{|\varphi|}{(1+|\psi|^2)^{1/2}}\in L_2(\RR)$. Then for any $\tau>0$ and any function $x\in L_2(\RR)$ such that $\psi{\cal F}[x]\in L_2(\RR)$ we have sharp inequalities (\ref{taik_gen_1}) and (\ref{taik_gen_2}).
\end{thm}

For the best approximation of the functional $\varphi(A)x(0)$ on the class $W^r_{2,2}(\RR)$ we have the following. If
\[
N=\left ( \displaystyle \displaystyle \int\limits_{\RR}\frac{|\varphi(s)|^2}{(1+\tau|\psi(s)|^2)^2}ds \right )^{1/2},
\]
 then
 \[
 E(N)=\tau \left ( \displaystyle \displaystyle \int\limits_{\RR}\frac{|\varphi(s)\psi(s)|^2}{(1+\tau|\psi(s)|^2)^2}ds \right )^{1/2},
 \]
and the extremal functional is
 \[
 g_\tau (x)=\displaystyle \int\limits_{\RR}\frac{\varphi(s)}{1+\tau |\psi(s)|^2}{\cal F}[x](s)ds.
 \]

\subsection{Results for functions of the differential operator in $L_2(\TT)$.}

Next let us take operator $A=\frac 1i\frac{d}{dt}$ in $L_2(\TT)$. The corresponding partition of unity is
\[
E(\beta)x(t)=\frac{1}{\sqrt{2\pi}}\displaystyle \sum_{n\in \ZZ\cap \beta}\hat x(n)e^{int},\qquad \beta\in {\cal B},
\]
where $\hat{x}(n)=\frac{1}{\sqrt{2\pi}}\displaystyle \int\limits_{-\pi}^\pi x(s)e^{-ins}ds$.
For the function $\varphi(A)$ of the operator $A$ we can also write
\[
\varphi(A)x(t)=\frac 1{\sqrt{2\pi}}\sum\limits_{n=-\infty}^{\infty}\varphi(n)\hat{x}(n)e^{int}.
\]

As usual, by $l_2(\ZZ)$ we denote the space of complex-valued sequences $\{ x_n\}=\{ x_n\}_{n\in \ZZ}$ such that $\sum_{n\in\ZZ}|x_n|^2<\infty$.
Domain of the operator $\varphi(A)$ consists of functions $x$ from $L_2(\TT)$ such that $\left\{ |\varphi(n)\hat{x}(n)|\right\}\in l_2(\ZZ)$.  If $x\in D(\psi(A))$ and condition $\left\{\frac{|\varphi(n)|}{(1+|\psi(n)|^2)^{1/2}}\right\}\in l_2(\ZZ)$ is satisfied, then it is easy to see that $\left\{\varphi(n) \hat{x}(n)\right\}\in l_1(\ZZ)$ and, hence, function $\varphi(A)x(t)$ is continuous.
Below we assume that these conditions are satisfied.

Let
\[
D_m(t)=\frac 1{\sqrt{2\pi}}\sum\limits _{n=-m}^{m}e^{int},\;\; m\in\mathbb{N},
\]
be the Dirichlet kernel. Then for any function $x\in L_{2}(\mathbb{T})$
\[
\frac{1}{2\pi}\displaystyle \int\limits _{0}^{2\pi}D_m(s-t)x(s)ds=S_{m}\left(x,t\right),
\]
where $S_{m}(x,t)$ is a Fourier partial sum for the function $x$ and, hence,
$$
(\varphi(A)x,D_m)=S_m(\varphi(A)x,0).
$$
Applying Theorem 2 to the functional on $L_2(\TT)$ defined by $D_m$, we obtain
\begin{equation}\label{201}
|(\varphi(A)x,D_m)|\leq N_{\varphi,\psi}(D_m;\tau) \|x\|_{L_2(\TT)}+\tau M{}_{\varphi,\psi}(D_m;\tau)\|\psi(A)x\|_{L_2(\TT)},
\end{equation}
Note that under the above conditions
$$
(\varphi(A)x,D_m)=S_m(\varphi(A)x,0)\to \varphi(A)x(0),\qquad m\to \infty ,
$$
$$
N_{\varphi,\psi}(D_m;\tau)\to\left\{\sum\limits_{n=-\infty}^\infty\frac{|\varphi(n)|^2}{(1+\tau |\psi(n)|^2)^2}\right\}^{1/2}
$$
and
$$
M_{\varphi,\psi}(D_m;\tau)\to\left\{\sum\limits_{n=-\infty}^\infty\frac{|\varphi(n)\psi(n)|^2}{(1+\tau |\psi(n)|^2)}\right\}^{1/2}.
$$
Taking $m\to \infty$ in (\ref{201}), we obtain
$$
|\varphi(A)x(0)|\le \left\{\sum\limits_{n=-\infty}^\infty\frac{|\varphi(n)|^2}{(1+\tau |\psi(n)|^2)^2}\right\}^{1/2} \| x\|_{L_2{(\TT)}}
$$
\begin{equation}\label{202}
\qquad\qquad\qquad+\tau \left\{\sum\limits_{n=-\infty}^\infty\frac{|\varphi(n)\psi(n)|^2}{(1+\tau |\psi(n)|^2)}\right\}^{1/2}\| \psi(A)x\|_{L_2(\TT)},
\end{equation}
and, therefore,
$$
\|\varphi(A)x\|_{L_\infty((\TT))}\le \left\{\sum\limits_{n=-\infty}^\infty\frac{|\varphi(n)|^2}{(1+\tau |\psi(n)|^2)^2}\right\}^{1/2} \cdot \| x\|_{L_2{(\TT)}}
$$
\begin{equation}\label{203}
\qquad\qquad\qquad+\tau \left\{\sum\limits_{n=-\infty}^\infty\frac{|\varphi(n)\psi(n)|^2}{(1+\tau |\psi(n)|^2)}\right\}^{1/2}\cdot  \| \psi(A)x\|_{L_2(\TT)}.
\end{equation}

Thus, we have proved the following theorem.
\begin{thm}\label{thm7}
Assume that functions $\varphi,\psi$ are continuous, satisfy conditions (\ref{condition}), and $\left\{\frac{|\varphi(n)|}{(1+|\psi(n)|^2)^{1/2}}\right\}\in l_2(\ZZ)$. Then for any $\tau>0$ for an arbitrary function $x\in L_2(\TT)$ such that $\left\{\psi(n)\hat{x}(n)\right\}\in l_2(\RR)$, sharp inequalities (\ref{202}) and (\ref{203}) hold.
\end{thm}

Theorem \ref{thm7} generalizes and sharpens mean-squared inequality of Shadrin \cite{Shadrin}.

For the best approximation of the functional $\varphi(A)x(0)$, the following is true.
If
\[
N=\left\{\sum\limits_{n=-\infty}^\infty\frac{|\varphi(n)|^2}{(1+\tau |\psi(n)|^2)^2}\right\}^{1/2},
\]
then
 \[
E(N)=\tau \left\{\sum\limits_{n=-\infty}^\infty\frac{|\varphi(n)\psi(n)|^2}{(1+\tau |\psi(n)|^2)}\right\}^{1/2},
 \]
and the extremal functional is
\[
g_\tau(x)=\sum\limits_{n=-\infty}^\infty\frac{\varphi(n)\hat{x}(n)}{1+\tau |\psi(n)|^2}.
\]

Note that the result about approximating this functional is new even in the case $\varphi(t)=t^k,\; \psi(t)=t^r,\;\; k,r\in\NN, \; k<r$.

\subsection{Orthogonal polynomials}

In this section we consider applications related to the partition of unity {connected with} classical orthogonal polynomials.

Let $I$ be the interval $(-1,1)$, real line $\RR$ or half line $\RR_+$. Let also for $\alpha,\beta>-1$ and $t\in I$
\[
h(t)=\left\{\begin{array}{clrc}
\mbox{$e^{-t^2}$} , & \mbox{ if $I=\RR$},\\
\mbox{$t^\alpha e^{-t}$} , & \mbox{ if $I=\RR_+$},\\
\mbox{$(1-t)^\alpha(1+t)^\beta$}, & \mbox{ if $I=(-1,1)$}.
\end{array}\right.
\]

By $L_{2,h}(a,b)$ we denote the Hilbert space of functions $x\in L_0(I)$ such that
\[
\left\| x\right\|_{L_{2,h}(a,b)}:=\left\{\displaystyle \int\limits_Ih(t)x^2(t)dt\right\}^{1/2}<\infty,
\]
with the inner product
\[
(x,y):=\displaystyle \int\limits_Ih(t)x(t)y(t)dt.
\]
By $\{F_n(x)\}$ we denote the system of orthogonal polynomials corresponding to the weight function $h(x)$, normalized so that $\| F_n\|_{L_{2,h}(I)}=1$.

As it is well-known (see, for instance, \cite[Ch. 2]{Suetin}), 
orthogonal polynomials $\{F_n(t)\}$ are solutions to the differential equation
\begin{equation}\label{difeq}
D(t)y''+(A(t)+D'(t))y'-\gamma_ny=0,
\end{equation}
where for $I=\RR$
\[
\begin{array}{clcr}
\mbox{$A(t)=-2t$},& \mbox{$D(t)\equiv 1$,} &\mbox{$\gamma_n=-2n$};
\end{array}
\]
for $I=\RR_+$
\[
\begin{array}{clcr}
\mbox{$A(t)=\alpha-t$,} &\mbox{ D(t)=t,}&\mbox {$\gamma_n=-n$};
\end{array}
\]
and for $I=(-1,1)$
\[
\begin{array}{clcr}
\mbox{$A(t)=\beta-\alpha-(\alpha+\beta)t$},&\mbox{$D(t)=1-t^2$,} &\mbox{$\gamma_n=-n(n+\alpha+\beta+1)$}.
\end{array}
\]

With the help of the sequence of orthogonal polynomials $\{F_n(t)\}$, we define the partition of unity in $L_{2,h}(I)$
\[
E(\beta)x(t)=\sum_{n\in \beta}x_nF_n(t),\qquad \beta\in {\cal B},
\]
where $x_n=(x,F_n)$. Then we have
\[
\displaystyle \int\limits_{\RR}dE(s)x(t)=\sum_{n=0}^\infty x_nF_n(t)=x(t).
\]

Differential equation (\ref{difeq}) can be re-written as follows
\[
 Ly=ny,
\]
where the differential operator $L$ is defined by
\[
Ly=D(t)y''+(A(t)+D'(t))y'-(\gamma_n-n)y.
\]
Then, on one hand, for any function $x\in D(L)$ we have
\[
Lx(t)=\sum_{n=0}^\infty nx_nF_n(t),
\]
 and on the other hand for the spectral integral  $
 \displaystyle \int_{\RR}tdE(t)$ and an arbitrary function $x\in D(L)$
we have
 \[
 \displaystyle \int\limits_{\RR}sdE(s)x(t)=\sum_{n=0}^\infty nx_nF_n(t).
 \]
 From here it follows that the defined partition of unity
generates operator $L$, and we can consider functions $\varphi(L)$  and $\psi(L)$ of this operator. So for the function $\varphi\in L_0(\RR)$
\begin{equation}\label{phiL}
 \varphi(L)x(t)=\displaystyle \int\limits_{\RR}\varphi(s)dE(s)x(t)=\sum_{n=0}^\infty \varphi(n)x_nF_n(t).
\end{equation}
Moreover, $x\in D(\varphi(L))$ if and only if $\{\varphi(n)x_n\}\in l_2(\ZZ_+)$.

Recall that if condition (\ref{condition}) is satisfied, then $D(\psi(L))\subset D(\varphi(L))$. Observe that if functions $\varphi,\psi$ are continuous, $x\in D(\psi(L))$, and the following holds
\begin{equation}\label{condition6}
\left\{\frac{|\varphi(n)|\| F_n\|_{C(I)}}{(1+|\psi(n)|^2)^{1/2}}\right\}\in l_2(\ZZ_+),
\end{equation}
then it is easy to see that
the series in the right-hand side of (\ref{phiL})
is uniformly convergent on $I$, and hence its sum $\varphi(L)x$ is continuous on $I$.

For $m\in \NN$ we set
\[
P_m(t,s)=\sum_{n=0}^mF_n(t)F_n(s).
\]
Clearly,
\[
(x(\cdot), P_m(t,\cdot))=\sum_{n=0}^m x_nF_n(t).
\]

Fixing $t\in I$, we define the functional $f_{t,m}$ on the space $L_{2,h}(I)$
\[
f_{t,m}(x)=(x(\cdot), P_m(t,\cdot)).
\]
Applying Theorem \ref{thm1}, we have
\begin{equation}\label{ort_pol_cons}
|f_{t,m}(\varphi(L)x)|\le N_{\varphi,\psi}(f_{t,m};\tau) \|x\|_{L_{2,h}(I)}+\tau M{}_{\varphi,\psi}(f_{t,m};\tau)\|\psi(A)x\|_{L_{2,h}(I)}
\end{equation}
Note that under the listed above conditions and with $m\to \infty $
$$
f_{t,m}(\varphi(L)x)\to \varphi(L)x(t),
$$
\begin{align*}
N_{\varphi,\psi}(f_{t,m};\tau)^2&=\displaystyle \int\limits_{\RR}\frac{|\varphi(s)|^2}{(1+\tau |\psi(s)|^2)^2}d(E(s)f_{t,m},f_{t,m})\\
&=\left(\displaystyle \int\limits_{\RR}\frac{|\varphi(s)|^2}{(1+\tau |\psi(s)|^2)^2}dE(s)f_{t,m},f_{t,m}\right)\\
&=\left(\sum\limits_{n=0}^\infty\frac{|\varphi(n)|^2}{(1+\tau |\psi(n)|^2)^2}(f_{t,m})_nF_n(t),f_{t,m}\right)\\
&=\sum\limits_{n=0}^m\frac{|\varphi(n)|^2}{(1+\tau |\psi(n)|^2)^2}F_n(t)^2\to \sum\limits_{n=0}^\infty\frac{|\varphi(n)|^2}{(1+\tau |\psi(n)|^2)^2}F_n(t)^2\\
\end{align*}
and
$$
M_{\varphi,\psi}(f_{t,m};\tau)^2\to\sum\limits_{n=0}^\infty\frac{|\varphi(n)\psi(n)|^2}{(1+\tau |\psi(n)|^2)^2}F_n(t)^2.
$$
Taking $m\to \infty$ in (\ref{ort_pol_cons}), we obtain
$$
|\varphi(L)x(t)|\le \left\{\sum\limits_{n=0}^\infty\frac{|\varphi(n)F_n(t)|^2}{(1+\tau |\psi(n)|^2)^2}\right\}^{1/2} \| x\|_{L_{2,h}{(I)}}
$$
\begin{equation}\label{1001}
\qquad\qquad\qquad+\tau \left\{\sum\limits_{n=0}^\infty\frac{|\varphi(n)\psi(n)F_n(t)|^2}{(1+\tau |\psi(n)|^2)^2}\right\}^{1/2}\| \psi(L)x\|_{L_{2,h}{(I)}}.
\end{equation}

Therefore, we have proved the following theorem.
\begin{thm}
Assume that continuous functions $\varphi,\psi$ satisfy conditions (\ref{condition}) and (\ref{condition6}). Then for any $\tau>0$, arbitrary function $x\in L_{2,h}{(I)}$ such that $\left\{\psi(n){x}(n)\right\}_{n=0}^\infty$ $\in l_2(\ZZ_+)$, sharp inequality (\ref{1001}) holds.
\end{thm}

For the best approximation of the functional $\varphi(L)x(t)$ the following is true. If
\[
N=\left\{\sum\limits_{n=0}^\infty\frac{|\varphi(n)F_n(t)|^2}{(1+\tau |\psi(n)|^2)^2}\right\}^{1/2},
\]
then
 \[
E(N)=\tau \left\{\sum\limits_{n=0}^\infty\frac{|\varphi(n)\psi(n)F_n(t)|^2}{(1+\tau |\psi(n)|^2)}\right\}^{1/2},
 \]
and the extremal functional is
\[
g_{t,\tau}(x)=\sum\limits_{n=0}^\infty\frac{\varphi(n){x}_n}{1+\tau |\psi(n)|^2}F_n(t).
\]


\end{document}